\documentclass{svmult}
\usepackage{amsmath}
\usepackage{amssymb}
\usepackage[english]{babel}
\usepackage{tikz}
\usetikzlibrary{intersections, shapes}
\usetikzlibrary{arrows}
\usepackage{subcaption}
\usepackage{booktabs,multirow}
\usepackage{bigdelim}
\usepackage[hidelinks]{hyperref}
\usepackage{array}
\usepackage{pgfplots}
\usepackage[leftcaption]{sidecap}
\sidecaptionvpos{figure}{t}
\usepackage{graphicx,bm,xcolor}

\begin{document}
\title*{IETI-DP for conforming multi-patch Isogeometric Analysis in three dimensions}
\titlerunning{}
\author{Rainer Schneckenleitner \inst{} and
	Stefan Takacs \inst{*}
}
\institute{
	Rainer Schneckenleitner \at RICAM, Austrian Academy of Sciences, Altenberger Straße 69, 4040 Linz, Austria
	\email{schneckenleitner@numa.uni-linz.ac.at}
	\and Stefan Takacs \at RICAM, Austrian Academy of Sciences, Altenberger Straße 69, 4040 Linz, Austria
	\email{stefan.takacs@ricam.oeaw.ac.at}\\
	\inst{*} Corresponding Author
}
\maketitle

\abstract{
We consider dual-primal isogeometric tearing and interconnection
(IETI-DP) solvers for multi-patch geometries in Isogeometric Analysis.
Recently, the authors have published a convergence analysis for those
solvers that is explicit in both the grid size and the spline degree
for conforming discretizations of two dimensional computational domains.
In the present paper, we shortly revisit these results and provide
numerical experiments that indicate that similar results may hold
for three dimensional domains. 
}

\section{Introduction}

We are interested in fast domain decomposition solvers for multi-patch
Isogoemtric Analysis (IgA; \cite{HughesCottrellBazilevs:2005}). We
focus on variants of FETI-DP solvers, see \cite{FarhatLesoinneLeTallecPiersonRixen:2001a,ToselliWidlund:2005a}
and references therein. Such methods have been adapted
to IgA in~\cite{KleissPechsteinJuttlerTomar:2012}, where the
individual patches of the multi-patch discretization are used
as subdomains for the solver. This method is sometimes referred to
as the dual-primal isogeometric tearing and interconnection (IETI-DP)
method. These methods are similar to Balancing Domain Decomposition
by Constraints (BDDC) methods, which have also been adapted for
IgA, see~\cite{BeiraoDaVeigaChoPavarinoScacci:2013} and follow-up papers by the same authors. The
similarity is outlined in~\cite{MandelDohrmannTezaur:2005a}.

Much progress for the IETI-DP methods has been made in the PhD-thesis
by C. Hofer, including the extension to various discontinuous Galerkin
formulations, see~\cite{Hofer:2016a}. Recently, the authors of the
paper at hand have extended the condition number bounds for the
preconditioned Schur complement system to be explicit not only
in the grid size but also in the chosen spline degree, see
\cite{SchneckenleitnerTakacs:2019} for the conforming case
and \cite{SchneckenleitnerTakacs:2020}
for an extension to the discontinuous Galerkin case. The analysis
follows the framework from~\cite{MandelDohrmannTezaur:2005a}. One key
ingredient for the analysis in~\cite{SchneckenleitnerTakacs:2019} has
been the construction of a bounded harmonic extension operator
for splines, which followed the ideas of~\cite{Nepomnyaschikh:1995}.
The analysis in~\cite{SchneckenleitnerTakacs:2019} treats the
two-dimensional case. As usual for FETI-like methods, the extension
of the analysis
to three dimensions is not effortless. The goal of this paper is
to demonstrate that the proposed method also performs well for
higher spline degrees in three dimensions.

The remainder of this paper is organized as follows. In Section~\ref{sec:2}, we introduce the model problem, discuss its discretization and
the proposed IETI-DP algorithm. In Section~\ref{sec:4}, numerical experiments for a three dimensional example are presented.

\section{Model problem and its solution}\label{sec:2}
We consider a standard Poisson model problem. Let $\Omega\subset \mathbb R^d$ be
a bounded Lipschitz domain. For given $f\in L_2(\Omega)$, we are
interested in solving for $u\in H^1(\Omega)$ such that
\[
		- \Delta u = f \quad \mbox{in}\quad \Omega
		\qquad \mbox{and}\qquad
		u=0\quad \mbox{on}\quad \partial\Omega
\]
holds in a weak sense. We assume that the closure of the 
computational domain $\Omega$ is the union of the closure of
$k$ non-overlapping patches $\Omega^{(k)}$ that are parametrized with
geometry functions
\[
			G_k : \widehat{\Omega} := (0,1)^d \rightarrow 
			\Omega^{(k)}:=G_k(\widehat{\Omega})
\]
such that for any $k\not=\ell$, the intersection $\overline{\Omega^{(k)}} \cap \overline{\Omega^{(\ell)}}$ is empty, a common vertex, a common edge, or
(in three dimensions) a common face (cf.~\cite[Ass.~2]{SchneckenleitnerTakacs:2019}). We assume that both, $\nabla G_k$
and $(\nabla G_k)^{-1}$, are in $L_\infty(\widehat{\Omega})$ for
all patches. For the analysis, we need a uniform bound on the
$L_\infty$-norm and a uniform bound on the number of neighbors of
each patch, cf.~\cite[Ass.~1 and 3]{SchneckenleitnerTakacs:2019}.

For each of the patches, we introduce a tensor B-spline discretization
on the parameter domain $\widehat\Omega$. The discretization is then
mapped to the physical patch $\Omega^{(k)}$ using the pull-back
principle. We use a standard basis as obtained by the Cox-de Boor
formula. We need a fully matching discretization, this means that
for each basis function that has a non-vanishing trace on one of
the interfaces, there is exactly one basis function on each of
the patches sharing this interface such that the traces of the basis
functions agree (cf.~\cite[Ass.~5]{SchneckenleitnerTakacs:2019}).
This is a standard assumption for any multi-patch setting that is not
treated using discontinuous Galerkin methods. For the analysis,
we assume quasi uniformity of grids within each
patch, cf.~\cite[Ass.~4]{SchneckenleitnerTakacs:2019}.

In the following, we explain how to the IETI-DP solver is set
up. Here, we loosely follow the notation used in the IETI-DP solution
framework that recently joined the public part of the G+Smo library.
We choose the patches as IETI subdomains.
We obtain patch-local stiffness matrices $A^{(k)}$ by evaluating
the bilinear forms $a^{(k)}(u,v) = \int_{\Omega^{(k)}} \nabla^\top u(x)\nabla v(x)
\mathrm dx$ with the basis functions for the corresponding patch.
We set up matrices $C^{(k)}$ such that their null space are the
coefficient vectors of the patch-local functions that vanish at
the primal degrees of freedom. In~\cite{SchneckenleitnerTakacs:2019},
we have considered corner values, edge averages, and 
the combination of both. In the three
dimensional case, we can choose corner values, edge averages,
face averages, and any combination thereof. We set
up fully redundant jump matrices $B^{(k)}$. We omit the
corner values if and only if the corners are chosen as primal
degrees of freedom. We setup the primal problem in the usual way,
i.e., we first, for $k=1,\ldots,K$, compute a basis by
\[
		\Psi^{(k)} := \begin{pmatrix} I&0 \end{pmatrix}
								(\widetilde{A}^{(k)})^{-1}
								  \begin{pmatrix} 0 \\ R_c^{(k)} \end{pmatrix},
		\quad \mbox{where}\quad
		\widetilde{A}^{(k)}
		:=
		\begin{pmatrix} A^{(k)} & (C^{(k)})^\top \\
								    C^{(k)} \end{pmatrix}
\]
and $R_c^{(k)}$ is a binary matrix that relates the primal constraints
(with their patch-local indices) to the degrees of freedom of
the primal problem (with their global indices) and set then
\[
		\widetilde A^{(K+1)} := \sum_{k=1}^K (\Psi^{(k)})^\top A^{(k)} \Psi^{(k)},
		\quad\mbox{and}\quad
		\widetilde B^{(K+1)} := \sum_{k=1}^K B^{(k)} \Psi^{(k)}.
\]
We consider the Schur complement problem
$
		F \underline \lambda = \underline g,
$
where
\[
		F := \sum_{k=1}^{K+1}
	\widetilde B^{(k)}			
		(\widetilde{A}^{(k)}) ^{-1}
	(\widetilde B^{(k)})^\top
	\quad\mbox{and}\quad
	\widetilde B^{(k)} :=
	\begin{pmatrix}
			B^{(k)} & 0
	\end{pmatrix}
	\mbox{ for } k=1,\ldots, K.
\]
The derivation of $\underline g$ is a patch-local preprocessing step.
We solve the Schur complement problem using a preconditioned
conjugate gradient (PCG) solver with the scaled Dirichlet
preconditioner
\[
		M_{\mathrm{sD}} := 
		\sum_{k=1}^K
	  B_{\Gamma}			
		D_k^{-1}
		\Big(
			A_{\Gamma\Gamma}^{(k)}
			- A_{\Gamma I}^{(k)} (A_{I I}^{(k)})^{-1} A_{I\Gamma}^{(k)}
		\Big)
		D_k^{-1}
	  ( B_{\Gamma} )^\top,
\]
where the index $\Gamma$ refers to the rows/columns of
$A^{(k)}$ and the columns of $B^{(k)}$ that refer to basis functions
with non-vanishing trace, the index $I$ refers to the remaining
rows/columns, and 
the matrix $D_k$ is a diagonal matrix defined based on
the principle of multiplicity scaling. For the analysis, it is
important that its coefficients are constant within each interface.
The solution $u$
itself is obtained from $\underline \lambda$ using the usual patch-local steps, cf.~\cite{SchneckenleitnerTakacs:2019}.

Under the presented assumptions, the condition number of the preconditioned Schur complement system
is in the two-dimensional case bounded by
\[
	C\, p \; 
	\left(1+\log p+\max_{k=1,\ldots,K} \log\frac{H_{k}}{h_{k}}\right)^2,
\]
where $p$ is the spline degree, $H_k$ is the
patch size, and $h_k$ the grid size,
see~\cite{SchneckenleitnerTakacs:2019}.

\section{Numerical results}\label{sec:4}

\begin{figure}[h]
	\centering
	\includegraphics[height=3.5cm]{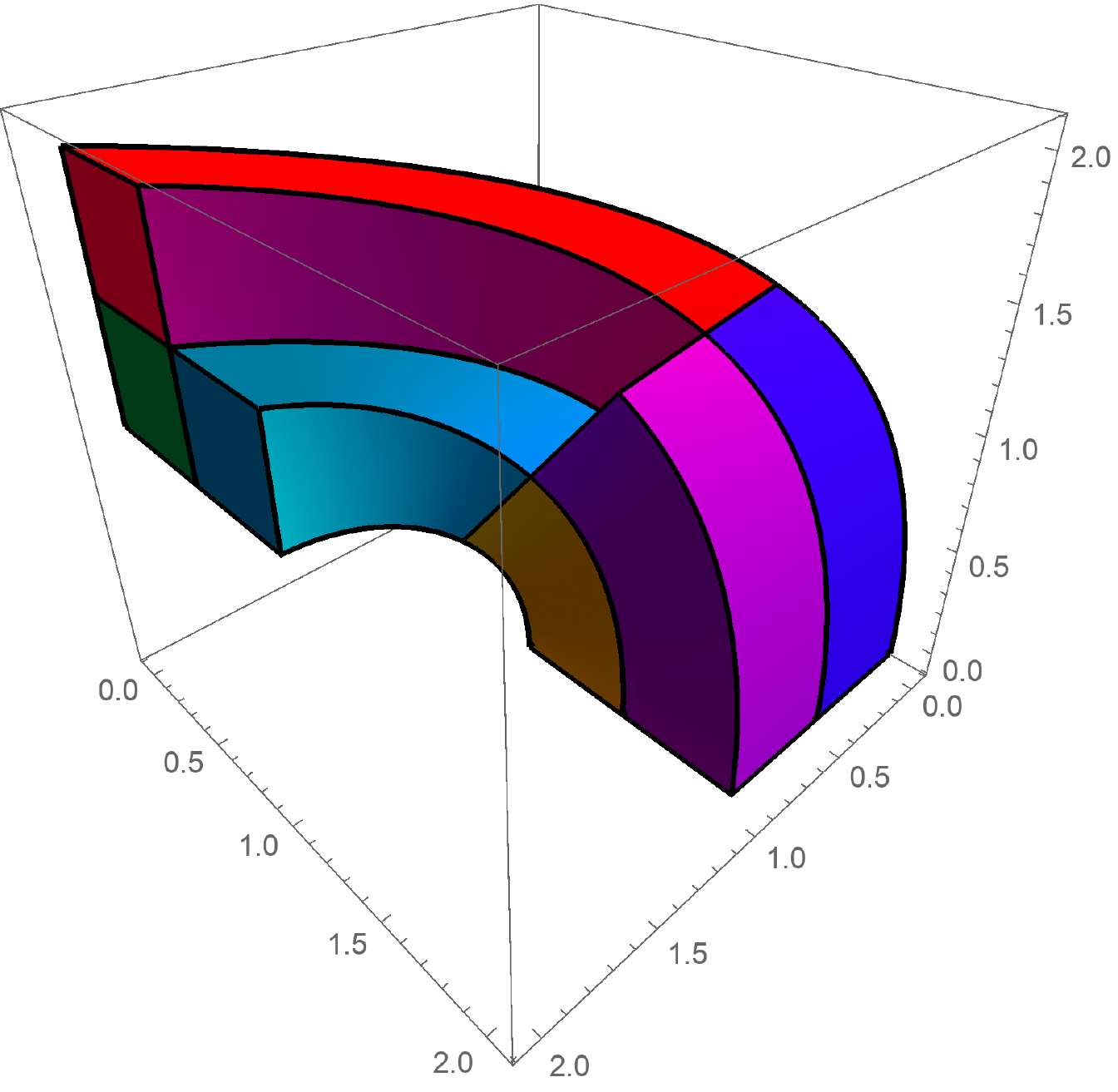}
	\caption{Computational domain}
	\label{fig:domain}
\end{figure}

In the following, we present numerical results for a three dimensional
domain and refer to the original paper~\cite{SchneckenleitnerTakacs:2019} for the two dimensional case.
The computational domain $\Omega$ is a twisted version of a Fichera
corner, see Fig.~\ref{fig:domain}.  The original geometry
consists of $7$ patches. We subdivide each patch uniformly into
$4\times 4\times 4$ patches to obtain a decomposition into
$448$~patches.

 We solve the model problem 
$- \Delta u(x,y,z)  = 3\pi^2 \sin(\pi x)\sin(\pi y)\sin(\pi z)$
for $(x,y,z)\in\Omega$
with homogeneous Dirichlet boundary conditions on $\partial\Omega$
by means of the IETI-DP solver outlined in the previous sections.
Within the patches, we consider tensor-product B-spline discretizations
of degree $p$ and maximum smoothness $C^{p-1}$. We consider several
grid sizes, the refinement level $r=0$ corresponds to a discretization
of each patch with polynomials. The next refinement levels $r=1,2,\ldots$
are obtained by uniform refinement.
All experiments have been carried out in the C++ library G+Smo\footnote{\url{https://github.com/gismo/gismo}, example file
\url{examples/ieti_example.cpp}.}
and have been executed on the Radon1 cluster\footnote{\url{https://www.ricam.oeaw.ac.at/hpc/}} in Linz. All computations have been performed with a single core.

\begin{figure}[ht]
\begin{center}
\resizebox{4.8cm}{4.8cm}{%
\begin{tikzpicture}
\begin{axis}[
xlabel={Refinement level $r$},
ylabel={Condition number $\kappa$},
xmin=1, xmax=4,
ymin=0, ymax=500,
xtick={1,2,3,4,5,6,7,8},
ytick={1,2,5,10,20,50,100,200,500},
legend pos=north west,
ymajorgrids=true,
grid style=dashed,
legend columns=3,
legend pos=north west,
legend image post style={only marks},
log ticks with fixed point,
xmode=log,
ymode=log
]

\addplot[color=blue,
mark=square]
coordinates {
(1,31.91)(2,79.95)(3,208.4)(4,509.37)
};
\addlegendentry{{\tiny V}};

\addplot[color=green,
mark=o]
coordinates {
(1,3.09)(2,4.64)(3,6.44)(4,8.61)
};
\addlegendentry{{\tiny E}};

\addplot[color=red,
mark=triangle]
coordinates {
(1,7.44)(2,10.69)(3,14.37)(4,18.37)
};
\addlegendentry{{\tiny F}};

\addplot[color=black,
mark=*]
coordinates {
(1,3.07)(2,4.62)(3,6.42)(4,8.59)
};
\addlegendentry{{\tiny V+E}};

\addplot[color=cyan,
mark=diamond]
coordinates {
(1,5.19)(2,8.56)(3,12.66)(4,17.40)
};
\addlegendentry{{\tiny V+F}};

\addplot[color=purple,
mark=otimes]
coordinates {
(1,3.06)(2,4.57)(3,6.36)(4,8.43)
};
\addlegendentry{{\tiny E+F}};

\addplot[color=orange,
mark=x]
coordinates {
(1,3.05)(2,4.55)(3,6.33)(4,8.42)
};
\addlegendentry{{\tiny V+E+F}};

\addplot[color=black,
dashed]
coordinates {
(1,10)(2,20)(3,30)(4,40)
};

\addplot[color=red,
dashed]
coordinates {
(1,40)(4,640)
};

\end{axis}
\end{tikzpicture}
}
\quad
\resizebox{4.8cm}{4.8cm}{%
\begin{tikzpicture}
\begin{axis}[
xlabel={Refinement level $r$},
ylabel={Solving times [sec]},
xmin=1, xmax=4,
ymin=0, ymax=5000,
xtick={1,2,3,4,5,6,7},
ytick={1,3.15,10,31.5,100,315,1000,3150},
legend pos=north west,
ymajorgrids=true,
grid style=dashed,
legend columns=3,
legend pos=north west,
legend image post style={only marks},
log ticks with fixed point,
ymode=log
]

\addplot[color=blue,
mark=square]
coordinates {
(1,7)(2,27)(3,200)(4,3873)
};
\addlegendentry{{\tiny V}};

\addplot[color=green,
mark=o]
coordinates {
(1,2.4)(2,7.2)(3,48)(4,679)
};
\addlegendentry{{\tiny E}};

\addplot[color=red,
mark=triangle]
coordinates {
(1,3.6)(2,11)(3,72)(4,1143)
};
\addlegendentry{{\tiny F}};

\addplot[color=black,
mark=*]
coordinates {
(1,2.4)(2,7.5)(3,47)(4,695)
};
\addlegendentry{{\tiny V+E}};

\addplot[color=cyan,
mark=diamond]
coordinates {
(1,3.2)(2,10)(3,69)(4,1015)
};
\addlegendentry{{\tiny V+F}};

\addplot[color=purple,
mark=otimes]
coordinates {
(1,2.5)(2,7.5)(3,55)(4,799)
};
\addlegendentry{{\tiny E+F}};

\addplot[color=orange,
mark=x]
coordinates {
(1,2.5)(2,7.6)(3,55)(4,804)
};
\addlegendentry{{\tiny V+E+F}};

\addplot[color=black,
dashed]
coordinates {
(1,10)
(2,127)
(3,3433)
(4,171447)
};

\end{axis}
\end{tikzpicture}
}
\end{center}
	\caption{Condition numbers and solving times for $p=3$
		\label{fig:1}}
\end{figure}
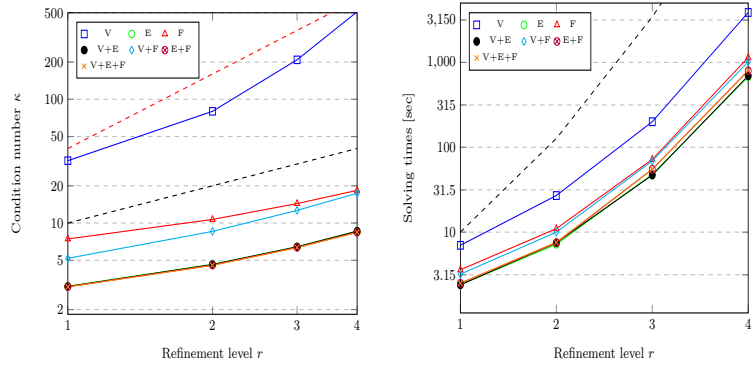

\begin{figure}[ht]
\begin{center}
\resizebox{4.8cm}{4.8cm}{%
\begin{tikzpicture}
\begin{axis}[
xlabel={Spline degree $p$},
ylabel={Condition number $\kappa$},
xmin=2, xmax=7,
ymin=0, ymax=500,
xtick={1,2,3,4,5,6,7,8},
ytick={1,2,5,10,20,50,100,200,500},
legend pos=north west,
ymajorgrids=true,
grid style=dashed,
legend columns=3,
legend pos=north west,
legend image post style={only marks},
log ticks with fixed point,
xmode=log,
ymode=log
]

\addplot[color=blue,
mark=square]
coordinates {
(2,41.77)(3,79.95)(4,122.43)(5,193.01)(6,226.67)(7,326.21)
};
\addlegendentry{{\tiny V}};

\addplot[color=green,
mark=o]
coordinates {
(2,3.86)(3,4.64)(4,5.33)(5,5.96)(6,6.69)(7,7.33)
};
\addlegendentry{{\tiny E}};

\addplot[color=red,
mark=triangle]
coordinates {
(2,9.49)(3,10.69)(4,11.77)(5,12.92)(6,14.16)(7,15.5)
};
\addlegendentry{{\tiny F}};

\addplot[color=black,
mark=*]
coordinates {
(2,3.84)(3,4.62)(4,5.31)(5,5.96)(6,6.69)(7,7.2)
};
\addlegendentry{{\tiny V+E}};

\addplot[color=cyan,
mark=diamond]
coordinates {
(2,6.8)(3,8.56)(4,10.02)(5,11.44)(6,12.86)(7,14.27)
};
\addlegendentry{{\tiny V+F}};

\addplot[color=purple,
mark=otimes]
coordinates {
(2,3.85)(3,4.57)(4,5.27)(5,5.94)(6,6.63)(7,7.41)
};
\addlegendentry{{\tiny E+F}};

\addplot[color=orange,
mark=x]
coordinates {
(2,3.83)(3,4.55)(4,5.26)(5,5.85)(6,6.53)(7,7.24)
};
\addlegendentry{{\tiny V+E+F}};

\addplot[color=red,
dashed]
coordinates {
(2,39.2)(7,480)
};

\addplot[color=black,
dashed]
coordinates {
(2,12)(7,42)
};

\end{axis}
\end{tikzpicture}
}
\quad
\resizebox{4.8cm}{4.8cm}{%
\begin{tikzpicture}
\begin{axis}[
xlabel={Spline degree $p$},
ylabel={Solving times [sec]},
xmin=2, xmax=7,
ymin=0, ymax=1000,
xtick={1,2,3,4,5,6,7},
ytick={1,2,5,10,20,50,100,200,500,1000},
legend pos=north west,
ymajorgrids=true,
grid style=dashed,
legend columns=3,
legend pos=north west,
legend image post style={only marks},
log ticks with fixed point,
xmode=log,
ymode=log
]

\addplot[color=blue,
mark=square]
coordinates {
(2,11)(3,27)(4,58)(5,130)(6,272)(7,548)
};
\addlegendentry{{\tiny V}};

\addplot[color=green,
mark=o]
coordinates {
(2,3.8)(3,7.2)(4,14)(5,28)(6,55)(7,104)
};
\addlegendentry{{\tiny E}};

\addplot[color=red,
mark=triangle]
coordinates {
(2,6.1)(3,11)(4,21)(5,40)(6,79)(7,150)
};
\addlegendentry{{\tiny F}};

\addplot[color=black,
mark=*]
coordinates {
(2,3.8)(3,7.5)(4,14)(5,28)(6,54)(7,101)
};
\addlegendentry{{\tiny V+E}};

\addplot[color=cyan,
mark=diamond]
coordinates {
(2,5.2)(3,10)(4,19)(5,39)(6,75)(7,143)
};
\addlegendentry{{\tiny V+F}};

\addplot[color=purple,
mark=otimes]
coordinates {
(2,4.1)(3,7.5)(4,15)(5,30)(6,57)(7,109)
};
\addlegendentry{{\tiny E+F}};

\addplot[color=orange,
mark=x]
coordinates {
(2,4.1)(3,7.6)(4,15)(5,28)(6,55)(7,105)
};
\addlegendentry{{\tiny V+E+F}};

\addplot[color=black,
dashed]
coordinates {
(2,11.8)(7,1184)
};

\end{axis}
\end{tikzpicture}
}
\end{center}
	\caption{Condition numbers and solving times for $r=2$
		\label{fig:2}}
\end{figure}
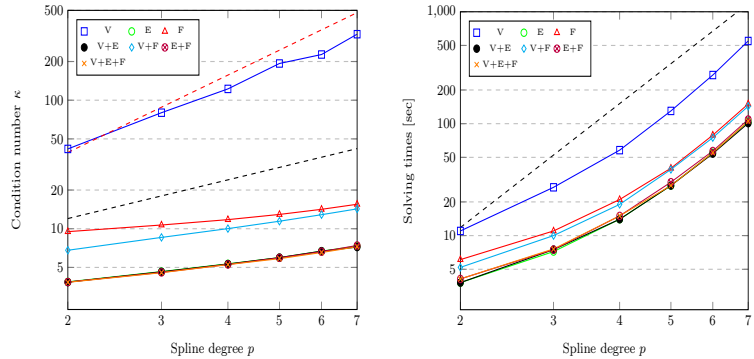


Concerning the choice of the primal degrees of freedom, we consider
all possibilities. For the two-dimensional case, the common choices
are the corner values, the edge averages, and a combination of both.
We have seen in~\cite{SchneckenleitnerTakacs:2019} that all approaches
work, typically the corner values better than the edge averages. As
expected, the combination of both yields the best results.
For the three dimensional case, we have more possibilities. We
report on these approaches in the
Tables~\ref{tab:V} (vertex values = V),~\ref{tab:E} (edge averages = E),~\ref{tab:F} (face averages = F),~\ref{tab:VE} (V+E),~\ref{tab:VF} (V+F),
and~\ref{tab:EF} (E+F).
The combination of all variants (V+E+F) is almost identical to the case V+E and only included in the
diagrams.
In any case, we report on the number of iterations (it) required
by the PCG solver to reduce the residual with a random
starting vector by a factor of $10^{-6}$ compared to the right-hand side. Moreover, we report on the condition numbers
$(\kappa)$ of the preconditioned system as estimated by the PCG solver.


\begin{table}[tb]
	\newcommand{\OoM}{\multicolumn{2}{c|}{OoM}}
	\newcommand{\OoMr}{\multicolumn{2}{c}{OoM}}
	\newcolumntype{L}[1]{>{\raggedleft\arraybackslash\hspace{-1em}}m{#1}}
	\centering
	\renewcommand{\arraystretch}{1.25}
		\begin{tabular}{L{1.5em}|L{2em}L{3em}|L{2em}L{3em}|L{2em}L{3em}|L{2em}L{3em}|L{2em}L{3em}|L{2em}L{3em}}
			\toprule
			& \multicolumn{2}{c|}{$p=2$} & \multicolumn{2}{c|}{$p=3$} 
			& \multicolumn{2}{c|}{$p=4$} & \multicolumn{2}{c|}{$p=5$}
			& \multicolumn{2}{c|}{$p=6$} & \multicolumn{2}{c }{$p=7$} \\
	$r$\;	 & it &$\kappa$\;& it &$\kappa$\;& it &$\kappa$\;& it &$\kappa$\;& it &$\kappa$\;& it &$\kappa$\;\\
			\midrule
  1\; & 33 &   14\; & 51 &   32\; & 64 &   45\; & 89 &   84\; &108 &  109\; &136 &  178\; \\
  2\; & 57 &   42\; & 79 &   80\; & 98 &  122\; &124 &  193\; &148 &  227\; &176 &  326\; \\
  3\; & 94 &  116\; &123 &  208\; &149 &  315\; &175 &  439\; &199 &  566\; &\OoMr        \\
  4\; &146 &  275\; &176 &  509\; &\OoM         &\OoM         &\OoM         &\OoMr        \\
   			\bottomrule
	\end{tabular}
	\caption{Iterations (it) and condition number ($\kappa$); Vertex (V)
		\label{tab:V}}
\end{table}

\begin{table}[tb]
	\newcommand{\OoM}{\multicolumn{2}{c|}{OoM}}
	\newcommand{\OoMr}{\multicolumn{2}{c}{OoM}}
	\newcolumntype{L}[1]{>{\raggedleft\arraybackslash\hspace{-1em}}m{#1}}
	\centering
	\renewcommand{\arraystretch}{1.25}
		\begin{tabular}{L{1.5em}|L{2em}L{3em}|L{2em}L{3em}|L{2em}L{3em}|L{2em}L{3em}|L{2em}L{3em}|L{2em}L{3em}}
			\toprule
			& \multicolumn{2}{c|}{$p=2$} & \multicolumn{2}{c|}{$p=3$} 
			& \multicolumn{2}{c|}{$p=4$} & \multicolumn{2}{c|}{$p=5$}
			& \multicolumn{2}{c|}{$p=6$} & \multicolumn{2}{c }{$p=7$} \\
	$r$\;	 & it &$\kappa$\;& it &$\kappa$\;& it &$\kappa$\;& it &$\kappa$\;& it &$\kappa$\;& it &$\kappa$\;\\
			\midrule
  1\; & 14 &  2.5\; & 17 &  3.1\; & 20 &  3.8\; & 23 &  4.4\; & 27 &  5.1\; & 29 &  5.5\; \\
  2\; & 18 &  3.9\; & 21 &  4.6\; & 23 &  5.3\; & 26 &  6.0\; & 29 &  6.7\; & 32 &  7.3\; \\
  3\; & 23 &  5.6\; & 25 &  6.4\; & 28 &  7.3\; & 30 &  8.0\; & 33 &  8.8\; &\OoMr        \\
  4\; & 27 &  7.5\; & 30 &  8.6\; &\OoM         &\OoM         &\OoM         &\OoMr        \\
			\bottomrule
	\end{tabular}
	\caption{Iterations (it) and condition number ($\kappa$); Edges (E)
		\label{tab:E}}
\end{table}

\begin{table}[tb]
	\newcommand{\OoM}{\multicolumn{2}{c|}{OoM}}
	\newcommand{\OoMr}{\multicolumn{2}{c}{OoM}}
	\newcolumntype{L}[1]{>{\raggedleft\arraybackslash\hspace{-1em}}m{#1}}
	\centering
	\renewcommand{\arraystretch}{1.25}
		\begin{tabular}{L{1.5em}|L{2em}L{3em}|L{2em}L{3em}|L{2em}L{3em}|L{2em}L{3em}|L{2em}L{3em}|L{2em}L{3em}}
			\toprule
			& \multicolumn{2}{c|}{$p=2$} & \multicolumn{2}{c|}{$p=3$} 
			& \multicolumn{2}{c|}{$p=4$} & \multicolumn{2}{c|}{$p=5$}
			& \multicolumn{2}{c|}{$p=6$} & \multicolumn{2}{c }{$p=7$} \\
	$r$\;	 & it &$\kappa$\;& it &$\kappa$\;& it &$\kappa$\;& it &$\kappa$\;& it &$\kappa$\;& it &$\kappa$\;\\
			\midrule
  1\; & 22 &  6.1\; & 26 &  7.4\; & 29 &  8.3\; & 33 &  9.5\; & 37 & 10.4\; & 41 & 11.5\; \\
  2\; & 29 &  9.5\; & 31 & 10.7\; & 34 & 11.8\; & 37 & 12.9\; & 42 & 14.2\; & 46 & 15.5\; \\
  3\; & 35 & 13.1\; & 38 & 14.4\; & 41 & 15.9\; & 43 & 17.0\; & 47 & 18.3\; &\OoMr        \\
  4\; & 41 & 17.1\; & 44 & 18.4\; &\OoM         &\OoM         &\OoM         &\OoMr        \\
			\bottomrule
	\end{tabular}
	\caption{Iterations (it) and condition number ($\kappa$); Faces (F)
		\label{tab:F}}
\end{table}

\begin{table}[tb]
	\newcommand{\OoM}{\multicolumn{2}{c|}{OoM}}
	\newcommand{\OoMr}{\multicolumn{2}{c}{OoM}}
	\newcolumntype{L}[1]{>{\raggedleft\arraybackslash\hspace{-1em}}m{#1}}
	\centering
	\renewcommand{\arraystretch}{1.25}
		\begin{tabular}{L{1.5em}|L{2em}L{3em}|L{2em}L{3em}|L{2em}L{3em}|L{2em}L{3em}|L{2em}L{3em}|L{2em}L{3em}}
			\toprule
			& \multicolumn{2}{c|}{$p=2$} & \multicolumn{2}{c|}{$p=3$} 
			& \multicolumn{2}{c|}{$p=4$} & \multicolumn{2}{c|}{$p=5$}
			& \multicolumn{2}{c|}{$p=6$} & \multicolumn{2}{c }{$p=7$} \\
	$r$\;	 & it &$\kappa$\;& it &$\kappa$\;& it &$\kappa$\;& it &$\kappa$\;& it &$\kappa$\;& it &$\kappa$\;\\
			\midrule
  1\; & 14 &  2.5\; & 17 &  3.1\; & 20 &  3.8\; & 22 &  4.3\; & 26 &  5.0\; & 28 &  5.4\; \\
  2\; & 18 &  3.8\; & 21 &  4.6\; & 23 &  5.3\; & 26 &  6.0\; & 29 &  6.7\; & 31 &  7.2\; \\
  3\; & 22 &  5.5\; & 25 &  6.4\; & 28 &  7.3\; & 30 &  8.0\; & 33 &  8.8\; &\OoMr        \\
  4\; & 27 &  7.5\; & 30 &  8.6\; &\OoM         &\OoM         &\OoM         &\OoMr        \\
			\bottomrule
	\end{tabular}
	\caption{Iterations (it) and condition number ($\kappa$); Vertices+Edges (V+E)
		\label{tab:VE}}
\end{table}

\begin{table}[tb]
	\newcommand{\OoM}{\multicolumn{2}{c|}{OoM}}
	\newcommand{\OoMr}{\multicolumn{2}{c}{OoM}}
	\newcolumntype{L}[1]{>{\raggedleft\arraybackslash\hspace{-1em}}m{#1}}
	\centering
	\renewcommand{\arraystretch}{1.25}
		\begin{tabular}{L{1.5em}|L{2em}L{3em}|L{2em}L{3em}|L{2em}L{3em}|L{2em}L{3em}|L{2em}L{3em}|L{2em}L{3em}}
			\toprule
			& \multicolumn{2}{c|}{$p=2$} & \multicolumn{2}{c|}{$p=3$} 
			& \multicolumn{2}{c|}{$p=4$} & \multicolumn{2}{c|}{$p=5$}
			& \multicolumn{2}{c|}{$p=6$} & \multicolumn{2}{c }{$p=7$} \\
	$r$\;	 & it &$\kappa$\;& it &$\kappa$\;& it &$\kappa$\;& it &$\kappa$\;& it &$\kappa$\;& it &$\kappa$\;\\
			\midrule
  1\; & 17 &  3.7\; & 22 &  5.2\; & 26 &  6.6\; & 30 &  7.8\; & 34 &  9.0\; & 38 & 10.1\; \\
  2\; & 25 &  6.8\; & 29 &  8.6\; & 32 & 10.0\; & 36 & 11.4\; & 40 & 12.9\; & 44 & 14.3\; \\
  3\; & 32 & 10.7\; & 36 & 12.7\; & 39 & 14.2\; & 42 & 15.6\; & 45 & 17.3\; &\OoMr        \\
  4\; & 39 & 15.2\; & 43 & 17.4\; &\OoM         &\OoM         &\OoM         &\OoMr        \\
			\bottomrule
	\end{tabular}
	\caption{Iterations (it) and condition number ($\kappa$); Vertices+Face (V+F)
		\label{tab:VF}}
\end{table}

\begin{table}[tb]
	\newcommand{\OoM}{\multicolumn{2}{c|}{OoM}}
	\newcommand{\OoMr}{\multicolumn{2}{c}{OoM}}
	\newcolumntype{L}[1]{>{\raggedleft\arraybackslash\hspace{-1em}}m{#1}}
	\centering
	\renewcommand{\arraystretch}{1.25}
		\begin{tabular}{L{1.5em}|L{2em}L{3em}|L{2em}L{3em}|L{2em}L{3em}|L{2em}L{3em}|L{2em}L{3em}|L{2em}L{3em}}
			\toprule
			& \multicolumn{2}{c|}{$p=2$} & \multicolumn{2}{c|}{$p=3$} 
			& \multicolumn{2}{c|}{$p=4$} & \multicolumn{2}{c|}{$p=5$}
			& \multicolumn{2}{c|}{$p=6$} & \multicolumn{2}{c }{$p=7$} \\
	$r$\;	 & it &$\kappa$\;& it &$\kappa$\;& it &$\kappa$\;& it &$\kappa$\;& it &$\kappa$\;& it &$\kappa$\;\\
			\midrule
  1\; & 14 &  2.5\; & 17 &  3.1\; & 20 &  3.8\; & 23 &  4.3\; & 27 &  5.0\; & 30 &  5.7\; \\
  2\; & 19 &  3.9\; & 21 &  4.6\; & 24 &  5.3\; & 27 &  5.9\; & 30 &  6.6\; & 33 &  7.4\; \\
  3\; & 23 &  5.5\; & 26 &  6.4\; & 29 &  7.2\; & 31 &  7.9\; & 33 &  8.5\; &\OoMr        \\
  4\; & 28 &  7.4\; & 31 &  8.4\; &\OoM         &\OoM         &\OoM         &\OoMr        \\
			\bottomrule
	\end{tabular}
	\caption{Iterations (it) and condition number ($\kappa$); Edges+Faces (E+F)
		\label{tab:EF}}
\end{table}

In Figure~\ref{fig:1}, the dependence on the refinement level is
depicted. Here, we have chosen the spline degree $p=3$ and have
considered all of the possibilities for primal degrees of freedom.
Here, we have 44\,965 
($r=1$), 133\,629 ($r=2$), 549\,037 ($r=3$), and
2\,934\,285 ($r=4$) degrees of freedom (dofs).
We observe that choosing only vertex values as primal degrees of
freedom leads to the largest condition numbers. We observe that
in this case the condition number grows like $r^2$ (the dashed red
line indicates the slope of such a growth). This corresponds to a growth like $(1+\log H/h)^2$, as
predicted by the theory for the two-dimensional case.
All other options yield
significantly better results, particularly those that include
edge averages. In these cases, the growth seems to be less
than linear in $r\approx \log H/h$ (the dashed black like shows
such a slope). In the right diagram, we can see that 
choosing a strategy with smaller condition numbers also yields a
faster method. Since the dimensions and
the bandwidths of the local stiffness matrices grow like
$(H_k/h_k)^3$ and $(H_k/h_k)^2$, respectively,
the complexity of the LU decompositions
grows like $\sum_{k=1}^K (H_k/h_k)^7$.
The complexity analysis indicates that they are the dominant factor.
The dashed black line indicates such a growth.

In Figure~\ref{fig:2}, the dependence on the spline degree is
presented, where we have chosen $r=2$. Here, the
number of dofs ranges from 66\,989 ($p=2$) to
549\,037 ($p=7$). Also in this picture,
we see that the vertex values perform worst and the edge averages
best. Again, we obtain a different asymptotic behavior
for the corner values. For those primal degrees of freedom,
the condition number grows like $p^2$ (the dashed red line indicates
the corresponding slope). All the other primal degrees of freedom
seem to lead to a growth that is smaller than linear in $p$
(the dashed black line indicates the slope of a linear growth).
Note that for the two-dimensional case, the theory predicts a
growth like $p (1+\log p)^2$. In the right diagram, we can see
that the solving times grow like $p^4$ (the dashed line shows the
corresponding slope). This seems to be realistic since the number
of non-zero entries of the stiffness matrix grows like $Np^d$, where
$N$ is the number of unknowns. For $d=3$, this yields in combination
with the condition number bound the observed rates.

Concluding, in this paper we have seen that the IETI method as
described in~\cite{SchneckenleitnerTakacs:2019} can indeed be
extended to the three dimensional case. As for finite elements,
only choosing vertex values is not enough.

\begin{acknowledgement}
	The first author was supported by the Austrian Science Fund (FWF): S117-03 and W1214-04. Also, the second author has received support from the Austrian Science	Fund (FWF): P31048.
\end{acknowledgement}

\bibliography{literature}

\end{document}